\documentclass[12pt]{amsart}
    \usepackage[utf8]{inputenc}
    
    \usepackage{amsmath,amssymb,amsfonts,amsthm}
	\usepackage{slashbox}

\author{Octavio A. Agustín-Aquino
}
\title{Some Patterns of Research in Mathematics}
\address{Instituto de Física y Matemáticas\\
Universidad Tecnológica de la Mixteca\\
Huajuapan de León, Oaxaca, México}
\date{May 21st, 2026}


\begin{document}

\maketitle

\begin{abstract}
We partially update Grossman's 2005 survey of patterns in mathematical research using a sample of 401 profiles from MathSciNet. The mathematical landscape has changed substantially: single-paper authors have reduced from $43$ \% to $32.42$ \%, collaboration has intensified, and the giant component of the co-authorship network has grown denser, with median Erd\H{o}s number dropping from $5$ to $4$. The data reveal the profession as markedly more collaborative, particularly among younger generations, although the tail of high productivity has also grown more extreme.
\end{abstract}

\section{Introduction}
I confess that the data I am going to present arise from the need to know where I stand with regard to the mass of the mathematicians. Jerrold Grossman \cite{Grossman2005} wrote a similar and more detailed survey back in 2005, for he had access to the full data of the MathSciNet database \cite{MathSciNet} but, alas, I do not. Further, even if I had a subscription, that does not mean it is available to explore it as extensively as he did. Yet I think an update is in order, even if it is only partial, for the panorama has changed substantially in two decades. Personally I draw three conclusions from the statistics of the present sample. First, the giant component of the mathematical collaboration graph appears to have become even more dominant, while the typical distances within it have been reduced. Second, comparisons with the whole MathSciNet population can be misleading without some attention to cohorts, since older and younger careers are mixed in the same archive. Third, collaboration now seems to occur earlier and more extensively, making the contemporary mathematical career structurally more of a collective endeavour than the one described by Grossman.



Thus first I probed the database. Since many MathSciNet author identifiers have six digits, I tested powers of ten and nearby values. Then I sampled uniformly from the interval of width $899000$. This yielded $343$ valid profiles from $500$ random identifiers, and continued sampling until obtaining $401$ valid ones.

With that sample, the estimated proportion of author profiles whose first publication predates 2005 is $64.1$ \%, and since in that year the total was about $3.34\times 10^{5}$, this yields an independent estimate of $5{.}21\times 10^{5}$ author profiles today.

\section{The statistics}

The main statistics of the sample appear in Table \ref{T:Statistics}. These were chosen because they are either central in Grossman's article or are reported by MathSciNet: number of publications, number of collaborators, Erdős-number data, MSC breadth, and unique citing authors.

The proportion of authors with just one paper is $32.42$ \%, with a $95$ \% confidence level interval of $[27.84,37.00]$. Hence it is quite certain that now it is less than the $43$ \% reported by Grossman. Similar assertions can be stated for the increase of the median number of papers, going from $2$ to $3$ (with a binomial sign test \cite{2023Sign} we have $z\approx -4{.}15$, $p\approx 1.6\times 10^{-5}$).

The average number of coauthors in the sample is $8{.}54$, which is much larger than the $2.94$ mean number of coauthors in $2005$; the subsample of those with at least one collaborator has an average of $9.76$ collaborators. The number of elements in the sample with at least one coauthor but with infinite Erd\H{o}s number is $36$. Hence approximately $9$ \% of the mathematicians of the sample are outside the giant component, and their mean number of coauthors is $1.75$. If I interpret correctly Grossman's data, the number of non-isolated mathematicians outside the giant component has decreased, for it was $13.35$ \%.

The number of unique citing authors (UCA) is not examined by Grossman, but it is quite telling that the maximum value in the sample is $2184$, while the median is $4$ and the mean is $63.51$. This reveals its distribution is highly skewed, with a very long right tail, which means UCA is better read as a rough indicator of the visibility structure of mathematical publications than as a description of a typical impact of a mathematician's work.

\begin{table}
\begin{center}
\begin{tabular}{|c|c|c|c|}
\hline
Number					& Mean 		& Std. Dev.	& Median\\
\hline
Papers 					& $12.47$	& $28.95$				& $3$\\
Coauthors				& $8.54$	& $18.55$				& $3$\\
UCA						& $63.51$	& $213.90$				& $4$\\
Erd\H{o}s				& $4.34$	& $0.99$				& $4$\\
\hline
\end{tabular}
\end{center}
\caption{Statistics from the sample. For the Erd\H{o}s number, it is calculated for those authors with a finite value.}
\label{T:Statistics}
\end{table}

It is also sobering to think about the deciles for the number of papers. Grossman reported the 6th, 7th, 8th and 9th deciles as $3$, $4$, $8$ and $17$. Twenty years later they have gone up to $4$, $6$, $13$ and $28.8$; they have increased at least $50$ \%. This seems to be due to the cumulative effect of longer careers, for if we restrict the calculations to a subsample of the mathematicians whose first paper was published in 2000 or afterwards, then the deciles drop to $3$, $6$, $10$ and $18.2$. It is very interesting that the deciles for coauthors increase in the lower regime (see Table \ref{T:Decileschicos}), suggesting that younger mathematicians tend to enter collaboration networks earlier than previous generations. This is accompanied by higher medium-range UCA values (the 4th, 5th and 6th deciles). The fact that their Erd\H{o}s numbers deciles remain the same further confirms the cohort's greater lean towards collaboration.

Regarding the ``top performance'' we have that the $95$ \% confidence interval\footnote{See \cite{Ialongo2019} to learn how the ranks for the interval are calculated.} for the 95th percentile for the number of papers within the whole sample is $[40,99]$, which clearly excludes the previous reported value of $30$; in the case of the subsample of those with first paper in 2000 or afterwards, the interval is $[21,45]$, which hints that the older value is still representative of the expected performance for younger mathematicians.

\begin{table}
\begin{center}
\begin{tabular}{|c|c|c|c|c|c|c|c|c|c|}
\hline
\backslashbox{Number}{Decile}	& 1		&2		&3		&4		&5		&6		&7		&8		&9\\
\hline
Papers 							& $1$	&$1$	&$1$	&$2$	&$3$	&$4$	&$6$	&$13$	&$28.8$\\
Coauthors						& $0$	&$1$	&$1$	&$2$	&$3$	&$4$	&$5$	&$9.3$	&$20$\\
UCA								& $0$	&$0$	&$0$	&$0$	&$4$	&$10$	&$19.2$	&$52.9$	&$131.8$\\
Erd\H{o}s						& $3$	&$4$	&$4$	&$4$	&$4$	&$5$	&$5$	&$5$	&$6$\\
\hline
\end{tabular}
\end{center}
\caption{Deciles from the sample. For the Erd\H{o}s number, it is calculated for those authors with a finite value.}
\label{T:Deciles}
\end{table}

\begin{table}
\begin{center}
\begin{tabular}{|c|c|c|c|c|c|c|c|c|c|}
\hline
\backslashbox{Number}{Decile}	& 1		&2		&3		&4		&5		&6		&7		&8		&9\\
\hline
Papers 							& $1$	&$1$	&$1$	&$2$	&$2$	&$3$	&$6$	&$10$	&$18.2$\\
Coauthors						& $1$	&$2$	&$2$	&$3$	&$3$	&$4.9$	&$6$	&$10$	&$15.6$\\
UCA								& $0$	&$0$	&$0$	&$3$	&$6$	&$10.9$	&$16$	&$29.4$	&$80$\\
Erd\H{o}s						& $3$	&$4$	&$4$	&$4$	&$4$	&$5$	&$5$	&$5$	&$6$\\
\hline
\end{tabular}
\end{center}
\caption{Deciles from the subsample ($n=199$) for mathematicians with first publication occurring in 2000 or afterwards. For the Erd\H{o}s number, it is calculated for those authors with a finite value.}
\label{T:Decileschicos}
\end{table}

The computation of the distribution of Erd\H{o}s numbers in the sample serves as another sanity check. Restricted to the giant component, it has mean $4.34$, standard deviation $0.99$ and median $4$, which is relatively close to the distribution reported by Grossman, who reported a mean of $4.69$, standard deviation $1.27$ and median $5$, for the MR collaboration graph. The small shift is plausibly explained by the densification of the network over the last two decades.

More precisely, if we restrict ourselves to the giant component ($n = 315$), $186$ authors have Erd\H{o}s number less than or equal to $4$ and $32$ have Erd\H{o}s number greater than $5$. A binomial sign test rejects the hypothesis that the median Erd\H{o}s number is $5$ in favor of a smaller value ($z \approx -10.4, p \approx 1.2\times 10^{-25}$). This provides strong evidence that the median Erd\H{o}s number has decreased to $4$ since the early 2000s. The reduction in the standard deviation also suggests that the distances have become more homogeneous, meaning the giant component has grown more connected.

The available data allow us to assign each mathematician a primary Mathematical Subject Classification (MSC) area, defined as the area in which they have published the largest number of papers, or, when this information is ambiguous, the area listed first in their record. For each individual, I recorded both the primary area and the total number of distinct MSC areas appearing in their publication record.

Among the 401 mathematicians considered, the most frequent primary areas (appearing in more than 1\% of profiles) were:

\begin{enumerate}
\item Operations research and mathematical programming (8.48\%)

\item Statistics (6.48\%)

\item Computer science (5.99\%)

\item Information and communication theory and circuits (5.24\%)

\item Biology and other natural sciences (4.74\%)

\item Numerical analysis (4.74\%)

\item Quantum theory (4.74\%)

\item Game theory, economics, finance, and social sciences (4.49\%)

\item Systems theory and control (4.24\%)

\item Combinatorics (3.74\%)

\item Fluid mechanics (3.49\%)

\item Partial differential equations (2.99\%)

\item Mechanics of deformable solids (2.74\%)

\item Optics and electromagnetic theory (2.74\%)

\item Probability theory and stochastic processes (2.74\%)

\item Ordinary differential equations (2.49\%)

\item Statistical mechanics and structure of matter (2.49\%)

\item Number theory (1.75\%)

\item Relativity and gravitational theory (1.75\%)

\item History and biography (1.25\%)

\item Operator theory (1.25\%)
\end{enumerate}

This list shows that identifying the mathematical profession with a small set of highly prestigious pure areas is misleading: a substantial part of the profiles lie in applied, statistical, computational, biological, physical, and operations-research areas.

The number of distinct MSC areas per individual provides a measure of disciplinary breadth. In the present sample, 50\% of mathematicians were associated with a single area, 70\% with at most two areas, 80\% with at most four areas, and 90\% with at most five areas. The maximum observed value was 18 areas. A comparison of the number of research areas of mathematicians with at least two publications is presented in Table \ref{T:Areas}. In the 2026 sample, there is a clear increase in disciplinary breadth, with a larger proportion of authors publishing in four or more areas.

\begin{table}
\begin{center}
\begin{tabular}{|c|c|c|}
\hline
Number of areas	& Grossman (2005), \% & Sample (2026), \%\\
\hline
1 & 38 &  35\\
2 & 32 &  24\\
3 & 14 & 10\\
4 & 7 & 10\\
5 & 4 & 7\\
6-10 & 5 & 9\\
11 or more & 0 & 4\\
\hline
\end{tabular}
\end{center}
\caption{Comparison of the number of areas in which authors with more than one paper published. The total of the sample does not equal 100 because some mathematicians were not assigned any MSC area.}
\label{T:Areas}
\end{table}

Thus, while most mathematicians remain highly specialized, a small fraction exhibit substantial interdisciplinary breadth. The resulting distribution displays a pronounced right tail, reflecting the presence of broadly active individuals who contribute to multiple research communities. This pattern is consistent with the observed structure of the coauthorship network and with the increasing prevalence of collaborative and interdisciplinary research.

\section{Conclusions}

It is interesting to observe, even through a relatively crude sampling method, the magnitude of the growth of the mathematical population. I learned a great deal of the human side of a mathematician's life simply by wandering through some profiles. For instance, I came across the profile of an émigré who landed at a public university and went on to publish over a hundred papers, with thousands of citations and some of the highest research distinctions. I was also amused to stumble upon the profiles of Helmut Hasse and Arnold Walfisz, while never picking one of a coauthor of Erd\H{o}s (as expected).

It is remarkable how extreme the tail behavior has become over the last twenty years: it is now substantially harder to reach the top decile in raw numbers. At the same time, younger generations appear to be doing no worse than their predecessors, but in a markedly more collaborative environment. The giant component is denser and more tightly interconnected than before. This increased interconnectedness promotes rapid integration into the network core and underscores the continuing importance of highly connected, cross-disciplinary mathematicians, in the tradition of Erd\H{o}s, in maintaining short collaboration paths.

As for the original motivation of this small project, it turns out that my own position in these distributions is better than I had expected, both in the full sample and within my cohort.

\bibliographystyle{amsplain}
\bibliography{notices_ams}

\end{document}